\documentclass[final]{siamltex}
\usepackage{amssymb,amsmath,graphicx}

 \title{THE ISENTROPIC EULER SYSTEM ADMITS
        SOME PLANE WAVE SUPERPOSITIONS}

 \author{ROBERT E. TERRELL\thanks{Mathematics Department, 
                   Cornell University, Ithaca, New York, 14853
                   ({\tt ret7@cornell.edu}).}}

\begin{document}

\maketitle

\begin{abstract}
 A class of differentiable solutions is proved for the 
 isentropic Euler equations in two and three space dimensions.
 The solutions are explicitly given in terms of 
 solutions to inviscid Burgers equations, and several 
 directions of propagation. The relative orientation of the 
 directions is critical. Within the directional constraints, 
 the Burgers solutions are arbitrary. 
 The several velocities add, and the pressures combine nonlinearly.
 These solutions cannot exist beyond the time 
 when shocks develop in any of the Burgers solutions.
\end{abstract}

\begin{keywords}
 Euler equations, Burgers equation, plane wave, isentropic, shock
\end{keywords}

\begin{AMS}
 35Q31, 76N15 
\end{AMS}

\pagestyle{myheadings}
\thispagestyle{plain}
\markboth{ROBERT E. TERRELL}{ISENTROPIC EULER SYSTEM}

\section{Main result}

Consider the isentropic Euler equations
$$
 u_t+u\cdot\nabla u +\rho^{-1}\nabla p = 0,
 \quad
 \rho_t+{\rm div}(\rho u) = 0,
 \quad
 p = k\rho^\gamma
$$
We assume $1<\gamma<3$, $k$ is constant, and set
$a=\frac{\gamma-1}{2}$.

\begin{theorem}
  Let 
 $v_j$ be unit vectors in 
 ${\mathbb R}^d$, $d=2$ or $3$, for which the dot products 
 $$
  v_i\cdot v_j = -a, \qquad i\ne j.
                                     \eqno{(*)}
 $$
 The number $N$ of such vectors is indicated in the table below.
 
 Further suppose that $f_j(s,t)$ are differentiable
 solutions to Burgers equation
 $$
   f_t+(1+a)ff_s = 0, \qquad s\in {\mathbb R}, \quad 0\le t < T
 $$
 Define
 $$
  u(x,t) = \sum_{j=1}^N f_j(x\cdot v_j,t)v_j,
  \qquad
  {\rm and\ }
  \rho = \Big(
       \frac{a}{\sqrt{k\gamma}}\sum_{j=1}^N f_j(x\cdot v_j,t)
         \Big)^{\frac{1}{a}}
 $$
 Then 
 $u$ and $\rho$ satisfy the isentropic Euler equations
 on this time interval and while $\sum f_j>0$. 

\end{theorem}
 
Note that the case $\gamma=2$ corresponds to the shallow water
model and you have
three vectors $v_k$ coplanar at $120$ degrees,  
while $\gamma=\frac{5}{3}$ corresponds to the monatomic gas
with four $v_k$ having the symmetry of a regular tetrahedron.

\begin{proof} 
We will write out the $d=3$ case, and the case $d=2$ can be obtained
by deleting the third component of all vectors.
Inspired by the treatment in Lax [3],
we
work with the symmetric hyperbolic form 
$$
 q_t+A_1q_{x_1}+A_2q_{x_2}+A_3q_{x_3}=0, 
 \qquad 
 q=\begin{bmatrix}u\\w\end{bmatrix}
$$
of the Euler equations
where
$q$ is a $4\times1$ vector
consisting of the velocities together with
$w=a^{-1}\sqrt{\gamma p/\rho}$, which is 
proportional to the sound speed.
That gives in the isentropic case
$\rho = (\frac{a}{\sqrt{k\gamma}} w)^{\frac{1}{a}}$.
The two forms of the Euler equations
 are equivalent for differentiable solutions
with $\rho>0$.
Here
$A_j= u_jI+awL_j$,
where $I$ is the $4\times 4$ identity matrix and
$$
 L_1 = \begin{bmatrix}0&0&0&1\\
                      0&0&0&0\\
                      0&0&0&0\\
                      1&0&0&0 
            \end{bmatrix},
 \quad 
 L_2 = \begin{bmatrix}0&0&0&0\\
                      0&0&0&1\\
                      0&0&0&0\\
                      0&1&0&0 
            \end{bmatrix},
 \quad
 L_3 = \begin{bmatrix}0&0&0&0\\
                      0&0&0&0\\
                      0&0&0&1\\
                      0&0&1&0 
            \end{bmatrix}
$$

We abbreviate
$\partial f_j/\partial s$ evaluated at $(x\cdot v_j,t)$ by $f_{js}$.
Component $i$ of vector $v_j$ is written $v_{ji}$. 
Also abbreviate 
$\partial f_j/\partial t(x\cdot v_j,t)$ by $f_{jt}$,
and $f_j(x\cdot v_j,t)$ by $f_j$.
Sums are from 1 to $N=2$, 3, or 4, depending on the number of vectors $v_k$.

We will need to know the eigenvectors of linear combinations of 
the $L_j$. These eigenvectors may be read from the calculation
$$
 \begin{bmatrix}0 & 0 & 0 & h \\
                0 & 0 & 0 & k \\
                0 & 0 & 0 & m \\
                h & k & m & 0 \end{bmatrix}
 \begin{bmatrix} h \\
                 k \\
                 m \\
               \pm1\end{bmatrix}
 =
 \pm \begin{bmatrix} h \\
                 k \\
                 m \\
               \pm1\end{bmatrix},
 \qquad
 \begin{bmatrix}0 & 0 & 0 & h \\
                0 & 0 & 0 & k \\
                0 & 0 & 0 & m \\
                h & k & m & 0 \end{bmatrix}
 \begin{bmatrix} h_0 \\
                 k_0 \\
                 m_0 \\
                 0   \end{bmatrix}
 = 0
$$
whenever $h^2+k^2+m^2 = 1$ and $hh_0+kk_0+mm_0 = 0$.

Now look for solutions of the form
$q(x,t) = \sum_k f_k(x\cdot v_k, t)z_k$
where 
constant vectors $v_k\in {\mathbb R}^3$ and 
$z_k\in {\mathbb R}^4$ 
are to be found.
Then
$$
 q_t+\sum_j (u_jI+awL_j)q_{x_j}
 =
 \sum_k \Big(
     f_{kt}+\sum_j (u_jI+awL_j)f_{ks}v_{kj}
        \Big)z_k
 $$ $$
 = \sum_k\Big(
         f_{kt}+f_{ks}u\cdot v_k +aw f_{ks}\sum_j(v_{kj}L_j)
         \Big) z_k
 $$
Now suppose we are looking for eigenvectors
 $\sum_j(v_{kj}L_j) z_k = \lambda_k z_k$.
As displayed above, we may either choose
$z_k = \begin{bmatrix}v_k\\\lambda_k\end{bmatrix}$
with 
$\lambda_k = \pm1$, or if
$\lambda_k = 0$ then the first three components of 
$z_k$ must be orthogonal to $v_k$.

With any such choices of eigenvectors then
$$
 q_t+\sum_j A_jq_{x_j}
 = 
 \sum_k\Big(
  f_{kt}+ f_{ks}\cdot(u\cdot v_k+aw\lambda_k)
     \Big) z_k
 $$ $$
 =
 \sum_k\Big(
  f_{kt}+ f_{ks}\cdot(
                 q\cdot\begin{bmatrix}v_k\\ a\lambda_k\end{bmatrix}
                     )
     \Big) z_k
 =
 \sum_k\Big(
  f_{kt}+ f_{ks}\cdot\Big(
         \sum_m f_m z_m\cdot\begin{bmatrix}v_k\\ a\lambda_k\end{bmatrix}
                     \Big)
 \Big) z_k
$$
We choose to make the dot products
$z_m\cdot\begin{bmatrix}v_k\\ a\lambda_k\end{bmatrix}=0$
for $k\ne m$,
which
decouples the system into the equations
$$
 f_{kt}+ 
   \Big(
 z_k\cdot\begin{bmatrix}v_k\\ a\lambda_k\end{bmatrix}
   \Big)
     f_k f_{ks}
 =0.
$$
If $\lambda_k=\pm1$
we have 
$z_k\cdot\begin{bmatrix}v_k\\ a\lambda_k\end{bmatrix}$
$=v_k\cdot v_k+a\lambda_k^2$
$=1+a$.
If
$\lambda_k=0$
then 
$z_k=\begin{bmatrix}v_k^\perp\\ 0\end{bmatrix}$ 
where $v_k^\perp$ is some vector perpendicular to $v_k$,
and $z_k\cdot\begin{bmatrix}v_k^\perp\\ 0\end{bmatrix}=0$,
so we need $f_k$ independent of $t$, as well as
$v_m^\perp\cdot v_k = 0$ for $k\ne m$.

Now we analyze the several cases of dot products and eigenvalues.
In the cases where some $\lambda_k = -1$, we 
replace $v_k$ by $-v_k$ and $f_k(s,t)$ by 
$-f_k(-s,t)$. This effectively replaces $-1$ by $+1$, and we
can assume from now on that all $\lambda_k\ge0$.

The most important case, and the one stated in the Theorem,
 is when all eigenvalues are $+1$.  
Let the $v_k$ be $N$ unit vectors with
all $v_k\cdot v_m = -a$ for  $k\ne m$, and all
$\lambda_k = 1$.  The $N$ is given in the table.
The decoupled equations for the $f_k$ are 
the inviscid Burgers [1] equation
$f_{kt}+ (1+a)f_kf_{ks}=0$.
The solutions are of the form
$q = \sum_{k=1}^N f_k(x\cdot v_k,t)\begin{bmatrix}v_k\\1\end{bmatrix}$.
This completes the proof. 
\qquad\end{proof}

Another possibility  is
that some eigenvalue is $0$.
Corresponding to each $0$ eigenvalue
 you may replace the term
$f_k(x\cdot v_k,t)\begin{bmatrix}v_k\\1\end{bmatrix}$
by
$g_k(x\cdot v_k)\begin{bmatrix}v_k^\perp\\0\end{bmatrix}$
where $g_k$ is any differentiable function,
provided that $v_k^\perp$ is perpendicular to $v_k$
and {\it all} the other $v_m$.

\begin{figure}[ht]
 \centerline{
   \begin{tabular}{|c|c|c|c|c|c|}
      \hline
 ${\mathbb R}^d$ & $1<\gamma<\frac{5}{3}$ & $\gamma=\frac{5}{3}$ & $\frac{5}{3}<\gamma<2$ & $\gamma=2$ & $2<\gamma<3$ \\
      \hline
    ${\mathbb R}^2$ & 2 &                   2 &         2 & 3 & 2 \\
    ${\mathbb R}^3$ & 3 &                   4 &         3 & 3 & 2 \\
      \hline
   \end{tabular} 
 }
\caption{
  The table shows the number $N$ of vectors $v_k$ 
 available for various $d$ and $\gamma$. 
 }
\label{fig0thetable}
\end{figure}

{\bf Remark on the time of existence.} Such configurations cannot generally
live beyond the time
when shocks develop in any of the $f_k$. 
For example, suppose a shock of speed $\sigma$
develops in $f_1$, and assume $\gamma = 1.4$. 
The jump condition on
density is
$[\rho]\sigma = [\rho u]\cdot v_1$
or
$$
 \left[\Big(
     \frac{a}{\sqrt{k\gamma}}\sum_k f_k
       \Big)^{\frac{1}{a}}
 \right]\sigma
 =
 \left[
       \Big(
     \frac{a}{\sqrt{k\gamma}}\sum_k f_k
       \Big)^{\frac{1}{a}}
   (f_1-af_2-af_3)
 \right].
$$
But this is not possible. Consider a line segment lying in
a plane level set of $f_3$ and within the shock plane. Along
this segment, $f_2$ will in general take a continuous range of
values, while $f_3$ is constant and $f_1$ has different
one-sided limits depending on the side of the shock plane from
which the segment is approached. Since $1/a=5$,
the jump condition is a polynomial identity in the 
values of $f_2$. This contradicts the fundamental theorem of algebra. 

Preliminary calculations done using
{\tt clawpack} \cite{clawpack}, \cite{LeVeque}
suggest that there is a distinction in the appearance of
pressure contours in two cases of crossing 
wave fronts shortly after breaking occurs, depending on whether the
angles between the fronts match equation $(*)$.

\end{document}